\documentclass[12pt]{article}
\usepackage{amsfonts}
\usepackage{amsmath,amssymb,amsfonts}

\newtheorem{theorem}{Theorem}

\newtheorem{lemma}[theorem]{Lemma}
\newtheorem{proposition}[theorem]{Proposition}

\begin{document}

\title{Centro-affine invariants and the canonical Lorentz metric on the
space of centered ellipses}
\author{Marcos Salvai\thanks{%
This work was partially supported by Conicet (PIP 112-2011-01-00670), Foncyt
(PICT 2010 cat 1 proyecto 1716) and Secyt Univ.\thinspace Nac.\thinspace C%
\'{o}rdoba.}}
\date{ }
\maketitle

\begin{abstract}
We consider smooth plane curves which are convex with respect to the origin.
We describe centro-affine invariants (that is, $GL_{+}\left( 2,\mathbb{R}%
\right) $-invariants), such as centro-affine curvature and arc length, in
terms of the canonical Lorentz structure on the three dimensional space of
all the ellipses centered at zero, by means of null curves of osculating
ellipses. This is the centro-affine analogue of the approach to conformal
invariants of curves in the sphere introduced by Langevin and O'Hara, using
the canonical pseudo Riemannian metric on the space of circles.
\end{abstract}

\footnote{%
 }

\smallskip

\noindent Runnning title: Centro-affine invariants via null curves of
ellipses

\section{Introduction}

Some years ago, R\'{e}mi Langevin and Jun O'Hara presented in \cite{LO'H} a
new approach to the classical subject of conformal length of curves in the
sphere, in terms of the canonical pseudo Riemannian structure on the space $%
\mathcal{C}$ of oriented circles: The circles osculating a curve $\alpha $
in the sphere define a null curve in $\mathcal{C}$, whose $\frac{1}{2}$%
-dimensional length provides, generically, a conformally invariant
parametrization of $\alpha $ (for the two-sphere $\mathcal{C}$ is Lorentz,
while it has signature $\left( 4,2\right) $ for the three-sphere). See also 
\cite{Fuster}. This line of thought can be traced back to Lie, Darboux and
Klein (see \cite{GLW}) and continued with the interpretation given by Robert
Bryant in \cite{bryant} of the standard conformally invariant 2-form on a
surface $M$ in $\mathbb{R}^{3}$ as the area of the surface in 5-dimensional
Lorentz space, determined by a certain family of tangent spheres to $M$.

In this paper we deal with an analogous situation: We consider smooth plane
curves which are convex with respect to the origin. The corresponding curves
of osculating ellipses centered at zero turn out to be null curves in $%
\mathcal{E}$, the three dimensional space of all ellipses centered at zero,
provided that $\mathcal{E}$ is endowed with a canonical Lorentz structure.
We use this notion to describe centro-affine, i.e.$\ GL_{+}\left( 2,\mathbb{R%
}\right) $-invariants.

I\ would like to thank the referee for reading the paper very carefully and
finding several typos.

\section{Centro-affine invariants of $0$-convex plane curves}

By a \emph{path }in the plane we understand an oriented embedded submanifold
included in $\mathbb{R}^{2}$ which is diffeomorphic to $\mathbb{R}$. Given a
path $c$ in the plane, an embedding $\alpha :I\rightarrow \mathbb{R}^{2}$
defined on the open interval $I$ with image $c$, such that $\alpha ^{\prime
} $ is positive with respect to the orientation of $c$, is called a \emph{%
parametrization} of $c$.

For $u,v\in \mathbb{R}^{2}$, denote $u\wedge v=\det \left( u,v\right) $. All
maps are supposed to be of class $C^{3}$, which we call smooth.\ 

\smallskip

A curve $\alpha :I\rightarrow \mathbb{R}^{2}$ is said to be \emph{convex
with respect to} $0$ (or briefly, $0$-\emph{convex}) if $\alpha \wedge
\alpha ^{\prime }$ and $\alpha ^{\prime }\wedge \alpha ^{\prime \prime }$
are both positive functions. In particular, $\alpha \left( t\right) \neq 0$
for all $t\in I$, $\alpha $ is regular, that is, $\alpha ^{\prime }$ never
vanishes, and $\alpha $ is traversed counterclockwise (we made this choice
for the sake of simplicity). A path $c$ in $\mathbb{R}^{2}$ is said to be $0$%
-convex if some (or equivalently, any) of its parametrizations is $0$-convex.

\smallskip

For instance, with a convenient orientation, a spiral centered at zero is $0$%
-convex, as well as an arc of the border of a strictly convex subset of the
plane containing the origin.

Let $\mathcal{P}$ be the set of paths in $\mathbb{R}^{2}$ which are $0$%
-convex. This set is invariant by the canonical action of the group $%
G:=GL_{+}\left( 2,\mathbb{R}\right) $ of linear isomorphisms of the plane
with positive determinant.

Any path $c\in \mathcal{P}$ admits a \emph{standard centro-affine
parame\-tri\-zation} $\alpha $, that is, a parame\-tri\-zation $\alpha $
such that 
\begin{equation}
\alpha ^{\prime \prime }=-\alpha +\tfrac{1}{2}\varkappa \alpha ^{\prime }
\label{cacurvature}
\end{equation}%
for some smooth function $\varkappa :I\rightarrow \mathbb{R}$, called the 
\emph{centro-affine curvature of} $\alpha $, which induces as usual a well
defined notion of centro-affine curvature on $c$. See \cite{Olver, Wilkens},
where the centro-affine curvature coincides up to a multiple with the one
given here and also 0-concave curves are considered simultaneously.

For the sake of completeness, we include the computation giving rise to $%
\varkappa $. Let $\beta $ a 0-convex curve and let $\alpha =\beta \left(
t\right) $ with $t^{\prime }>0$. Then $\alpha ^{\prime }=\beta ^{\prime
}\left( t\right) t^{\prime }$, $\alpha ^{\prime \prime }=\beta ^{\prime
\prime }\left( t\right) \left( t^{\prime }\right) ^{2}+\beta ^{\prime
}\left( t\right) t^{\prime \prime }$. We have $\beta ^{\prime \prime
}=a\beta +b\beta ^{\prime }$ for some functions $a,b$ with $a<0$. If the
function $t$ satisfies the equation $\left( t^{\prime }\right) ^{2}a\left(
t\right) =-1$, then (\ref{cacurvature}) holds with $\varkappa =2b\left(
t\right) t^{\prime }+2t^{\prime \prime }/t^{\prime }$ and so $\alpha $ is a
standard centro-affine reparametrization of $\beta $. It is an easy to
verify fact that a $0$-convex path is an arc of an ellipse if and only if $%
\varkappa \equiv 0$.

There is a broader notion of centro-affine arc length:\ Suppose $\mathcal{P}%
_{o}$ is a subset of $\mathcal{P}$ closed under the action of $G$. Any map
defined on $\mathcal{P}_{o}$ assigning to $c\in \mathcal{P}_{o}$ a nowhere
vanishing positively oriented $1$-form $\tau _{c}$ on $c$ is called a \emph{%
centro-affine arc length element on }$\mathcal{P}_{o}$, provided that $\tau
_{c}=g^{\ast }\tau _{gc}$ for any $g\in G$. This induces a $G$-invariant way
of measuring length of arcs of $0$-convex paths. In the analogous three
dimensional conformal setting, the conformal arc length element is not
defined for any path, but only for the so called vertex free paths, having
parametrizations with $\left( \kappa ^{\prime }\right) ^{2}+\left( \kappa
\tau \right) ^{2}>0$, where $\kappa $ and $\tau $ denote the curvature and
torsion (see Definition 1.2 and Theorem 7.3 in \cite{LO'H}). Theorem \ref%
{caelement} below involves this notion.

\section{The canonical Lorentz metric on the space of centered ellipses}

A subset $E$ of $\mathbb{R}^{2}$ is an ellipse centered at zero if there
exist an orthonormal basis $u,v$ of $\mathbb{R}^{2}$ and positive numbers $%
a,b$ such that 
\begin{equation}
E=\left\{ xu+yv\mid \frac{x^{2}}{a^{2}}+\frac{y^{2}}{b^{2}}=1\right\} \text{.%
}  \label{elipseuv}
\end{equation}%
We will consider only ellipses centered at zero, so in the following we
sometimes call them just ellipses.

Let $\mathcal{E}$ be the set of all ellipses in the plane centered at zero
(with axes not necessarily parallel to the coordinate axes) and let $%
\mathcal{S}_{+}$ be the manifold of all positive definite symmetric $2\times
2$ matrices. Among the several ways of identifying $\mathcal{E}$ with $%
\mathcal{S}_{+}$ we choose the following: 
\begin{equation}
F:\mathcal{S}_{+}\rightarrow \mathcal{E}\text{, \ \ \ \ \ \ \ }F\left(
A\right) =E_{A}\text{,}  \label{isometry}
\end{equation}%
where 
\begin{equation}
E_{A}=\left\{ z\in \mathbb{R}^{2}\mid \left\langle A^{-1}z,z\right\rangle
=1\right\} =A^{1/2}S^{1}=\left\{ A^{1/2}z\mid \left\vert z\right\vert
=1\right\} \text{,}  \label{ellipse}
\end{equation}%
since it is equivariant with respect to the canonical smooth transitive
actions of the group $G$ on $\mathcal{S}_{+}$ and $\mathcal{E}$, given by $%
g\cdot A=gAg^{T}$ and $g\cdot E=g\left( E\right) $, respectively (the
superscript $T$ denotes transpose). Notice that $E_{A}$ is equal to $E$ as
in (\ref{elipseuv}) provided that $Au=a^{2}u$ and $Av=b^{2}v$.

Now, $\mathcal{S}_{+}$ is an open set in the three dimensional vector space $%
\mathcal{S}$ of $2\times 2$ real symmetric matrices. We consider on $%
\mathcal{S}_{+}$ the unique $G$-invariant Lorentz structure on $\mathcal{S}%
_{+}$ whose norm at the identity $I$ is given by 
\begin{equation*}
\left\langle X,X\right\rangle =_{\text{def}}\left\Vert X\right\Vert =-\det X%
\text{, \ \ \ \ \ \ \ for }X\in T_{I}\mathcal{S}_{+}\cong \mathcal{S}\text{.}
\end{equation*}%
Equivalently, $\left\Vert \left( A,X\right) \right\Vert =-\det \left(
A^{-1}X\right) $ for any $\left( A,X\right) \in T\mathcal{S}_{+}\cong 
\mathcal{S}_{+}\times \mathcal{S}$.

We define the future pointing cone in $T_{I}\mathcal{S}_{+}$ as the set of
all $X\in \mathcal{S}$ with $\left\Vert X\right\Vert \leq 0$ such that
either $X_{11}$ or $X_{22}$ is positive. This induces a temporal orientation
on $\mathcal{S}_{+}$, which is invariant by the action of $G$.

A tangent vector $X$ of a Lorentz manifold is called spatial, temporal or
null (or light-like), if $\left\Vert X\right\Vert $ is positive, negative or
zero, respectively.

Consider on the space of ellipses $\mathcal{E}$ the Lorentz metric copied
from that in $\mathcal{S}_{+}$ above via the bijection (\ref{isometry}).

\begin{proposition}
The $G$-invariant metric on $\mathcal{E}$ defined above is isometric to $%
\mathcal{H}\times _{-}\mathbb{R}$, the warped product of the hyperbolic
plane $\mathcal{H}$ of constant curvature $-1$ with $\mathbb{R}$ with
warping function $\mathcal{H}\rightarrow \mathbb{R}$ constant and equal to $%
-1$.
\end{proposition}

\noindent \textbf{Proof. }Let $H=SL\left( 2,\mathbb{R}\right) =\left\{ A\in 
\mathbb{R}^{2\times 2}\mid \det A=1\right\} $ endowed with the bi-invariant
Lorentz metric defined at the identity by $\left\Vert X\right\Vert _{H}=%
\frac{1}{2}~$tr~$\left( X^{2}\right) $ (tr~$X=0$), which is a multiple of
the Killing form of $H$. It is well known that there is a pseudo Riemannian
submersion from $H$ onto $\mathcal{H}$ with isotropy group $SO\left(
2\right) $.

Now, $F:H\times \mathbb{R}\rightarrow G$ defined by $F\left( A,x\right)
=e^{x}A$ is a Lie group isomorphism satisfying%
\begin{equation*}
\left\Vert dF_{\left( I,0\right) }\left( X,x\left. \tfrac{d}{dt}\right\vert
_{0}\right) \right\Vert =-\det \left( X+xI\right) =\frac{1}{2}~\text{tr~}%
\left( X^{2}\right) -x^{2}=\left\Vert X\right\Vert _{H}-x^{2}\text{.}
\end{equation*}%
Hence $F$ is an isometry between $H\times _{-}\mathbb{R}$ and $G$.
Considering the quotient by $SO\left( 2\right) \times \left\{ 0\right\}
\simeq SO\left( 2\right) $, the proposition follows. \hfill $\square $

\bigskip

One can see $\mathcal{E}$ as the set of curves in the plane congruent to the
circle via the group $G$. Let $K$ be a Lie group acting on a manifold $N$.
Canonical $K$-invariant pseudo Riemannian metrics on spaces of $K$-congruent
curves in $N$ have proved to be useful in the study of foliations of $N$ by
such curves (see for instance \cite{SalvaiBLMS, GSmz,GSpams}).

\bigskip

Although the $G$-invariant metric on $\mathcal{S}_{+}$ is relevant for the
nature of the results, for some computations it will be convenient to
consider on $\mathcal{S}_{+}$ the constant Lorentz structure $g$ whose
associated norm is $\left\Vert X\right\Vert =-\det \left( X\right) $ (notice
that $\mathcal{S}_{+}$ is an open subset of the vector space $\mathcal{S}$).
The $G$-invariant metric $\bar{g}$ defined above is conformally equivalent
to $g$, and $\bar{g}=\phi ^{-2}g$, where $\phi :\mathcal{S}_{+}\rightarrow 
\mathbb{R}$ is given by $\phi \left( A\right) =\sqrt{\det \left( A\right) }$.

\begin{lemma}
\label{DbarD}Let $M$ be a smooth manifold and let $\bar{g}$ and $g$ be two
conformally equivalent pseudo Riemannian metrics on $M$. If $\gamma $ is a
smooth curve in $M$ which is null for $\bar{g}$ \emph{(}or equivalently, for 
$g$\emph{)}, then 
\begin{equation*}
\left\Vert \frac{\bar{D}\gamma ^{\prime }}{dt}\right\Vert =\left\Vert \frac{%
D\gamma ^{\prime }}{dt}\right\Vert \text{,}
\end{equation*}%
where $\frac{\bar{D}}{dt}$ and $\frac{D}{dt}$ denote the covariant
derivatives along $\gamma $ associated with $\bar{g}$ and $g$, respectively.
\end{lemma}

\noindent \textbf{Proof.} It suffices to show that $\left\Vert \bar{\nabla}%
_{X}X\right\Vert =\left\Vert \nabla _{X}X\right\Vert $ for any null local
vector field on $M$. Now, we have from the proof of Proposition 2.2 in \cite%
{Hans-Bert} that for any vector field on $M$, 
\begin{equation*}
\bar{\nabla}_{X}X=\nabla _{X}X-2X\left( \log \phi \right) X+\left\Vert
X\right\Vert \operatorname{grad}\left( \log \phi \right) \text{,}
\end{equation*}%
where $\bar{g}=\phi ^{-2}g$, with $\phi :M\rightarrow \mathbb{R}$. Hence, if 
$X$ is null, 
\begin{equation*}
\left\Vert \bar{\nabla}_{X}X\right\Vert =\left\Vert \nabla _{X}X\right\Vert
-4 X\left( \log \phi \right) \left\langle \nabla _{X}X,X\right\rangle
=\left\Vert \nabla _{X}X\right\Vert \text{,}
\end{equation*}%
as desired. The last equality follows from the fact that $\left\langle
\nabla _{X}X,X\right\rangle =0$ since $\left\Vert X\right\Vert $ is
constant. \hfill $\square $

\section{Null curves of osculating ellipses}

Given a regular curve $\alpha :I\rightarrow \mathbb{R}^{2}$, the (Euclidean)
curvature of $\alpha $ is the real function $\kappa $ on $I$ defined by $%
\kappa \left( t\right) =\left( \alpha ^{\prime }\wedge \alpha ^{\prime
\prime }\right) /\left\vert \alpha ^{\prime }\right\vert ^{3}$. This induces
a well defined notion of curvature on a path in $\mathbb{R}^{2}$.

For each $z$ on a $0$-convex path $c$ there exists exactly one ellipse $E$
osculating $c$ at $z$. This means the following:\ Suppose that $\alpha $ is
a parametrization of $c$ with $\alpha \left( t_{o}\right) =z$, then $E$ is
the ellipse having a counterclockwise parametrization $\varepsilon $ such
that $\varepsilon \left( 0\right) =z$, $\varepsilon ^{\prime }\left(
0\right) $ is a multiple of $\alpha ^{\prime }\left( t_{o}\right) $ and $%
\kappa _{\varepsilon }\left( 0\right) =\kappa _{\alpha }\left( t_{o}\right) $%
, where $\kappa _{\varepsilon }$ and $\kappa _{\alpha }$ are the curvatures
of $\varepsilon $ and $\alpha $, respectively.

\bigskip

The first assertion of the next theorem is the centro-affine analogue of
Theorems 5.1 and 7.2 in \cite{LO'H}, within the conformal setting, for $%
S^{2} $ and $S^{3}$, respectively.

\begin{theorem}
\label{TeoP}Let $\alpha :I\rightarrow \mathbb{R}^{2}$ be a parametrization
of a $0$-convex path $c$. For each $t\in I$, let $E\left( t\right) $ be the
ellipse centered at the origin osculating $\alpha $ at $t$. Then the curve $%
E:I\rightarrow \mathcal{E}$ is light-like.

If $\alpha $ is a standard centro-affine parametrization of $c$ and $%
\varkappa :I\rightarrow \mathbb{R}$ is the centro-affine curvature of $%
\alpha $, then%
\begin{equation}
\varkappa ^{2}\left( t\right) =\left\Vert \frac{D}{dt}E^{\prime }\left(
t\right) \right\Vert  \label{igualdad}
\end{equation}%
for all $t\in I$. Moreover, $\varkappa $ is positive or negative depending
on whether $E$ is future or past directed.
\end{theorem}

\noindent \textbf{Proof. }If\textbf{\ }$\alpha \circ s$ is a
reparametrization of $\alpha $, then $E\circ s$ is the curve of osculating
ellipses to $\alpha $, which is null if and only if $E$ is null. Hence, we
may suppose that $\alpha $ is a standard centro-affine parametrization, that
is, $\alpha $ satisfies (\ref{cacurvature}).

Given $t_{o}\in I$, we verify that $\left\Vert E^{\prime }\left(
t_{o}\right) \right\Vert =0$. We may assume additionally, without loss of
generality, that $t_{o}=0,$%
\begin{equation}
\alpha _{0}=e_{1}\text{,\ \ \ \ \ \ }\alpha _{0}^{\prime }=e_{2}
\label{normaliz}
\end{equation}%
(by considering $\bar{\alpha}\left( t\right) =A\alpha \left( t-t_{o}\right) $%
, where $A\in G$ satisfies $A\alpha _{t_{o}}=e_{1}$ and $A\alpha
_{t_{o}}^{\prime }=e_{2}$ and using the $G$-invariance of the statement).

Given $t$ in the domain of $\alpha $, we see first that 
\begin{equation*}
u\mapsto \varepsilon _{t}\left( u\right) =\left( \cos u\right) \alpha
_{t}+\left( \sin u\right) \alpha _{t}^{\prime }
\end{equation*}%
parametrizes $E\left( t\right) $. Clearly $\varepsilon _{t}$ is a
counterclockwise parametrization of an ellipse. We have that 
\begin{equation*}
\varepsilon _{t}^{\prime }\left( u\right) =-\left( \sin u\right) \alpha
_{t}+\left( \cos u\right) \alpha _{t}^{\prime }
\end{equation*}%
and so $\varepsilon _{t}^{\prime \prime }=-\varepsilon _{t}$. In particular, 
$\varepsilon _{t}\left( 0\right) =\alpha _{t}$, $\varepsilon _{t}^{\prime
}\left( 0\right) =\alpha _{t}^{\prime }$ and $\varepsilon _{t}^{\prime
\prime }\left( 0\right) =-\alpha _{t}$. So, for $E\left( t\right) $ to
osculate $\alpha $ at $t$ it suffices to check that $\kappa _{\varepsilon
_{t}}\left( 0\right) =\kappa _{\alpha }\left( t\right) $, where the left and
right hand side are the curvatures of $\varepsilon _{t}$ and $\alpha $ at $%
u=0$ and $t$, respectively. Indeed, 
\begin{equation*}
\kappa _{\alpha }\left( t\right) =\frac{\alpha _{t}^{\prime }\wedge \alpha
_{t}^{\prime \prime }}{\left\vert \alpha _{t}^{\prime }\right\vert ^{3}}=%
\frac{\alpha _{t}^{\prime }\wedge \left( -\alpha _{t}+\frac{1}{2}\varkappa
_{t}\alpha _{t}^{\prime }\right) }{\left\vert \alpha _{t}^{\prime
}\right\vert ^{3}}=\frac{\alpha _{t}\wedge \alpha _{t}^{\prime }}{\left\vert
\alpha _{t}^{\prime }\right\vert ^{3}}=\frac{\varepsilon _{t}^{\prime
}\left( 0\right) \wedge \varepsilon _{t}^{\prime \prime }\left( 0\right) }{%
\left\vert \varepsilon _{t}^{\prime }\left( 0\right) \right\vert ^{3}}%
=\kappa _{\varepsilon _{t}}\left( 0\right) \text{.}
\end{equation*}

Now we find the curve in $\mathcal{S}_{+}$ associated with the curve $%
E\left( t\right) $ in $\mathcal{E}$, according to the isometry (\ref%
{isometry}). Let $A_{t}=\left( \alpha _{t},\alpha _{t}^{\prime }\right) \in 
\mathbb{R}^{2\times 2}$, where $\alpha _{t}$ and $\alpha _{t}^{\prime }$ are
column vectors. We have that $E\left( t\right) =A_{t}S^{1}$. For simplicity,
in the following we omit writing $t$. The polar decomposition of $A$ is
given by%
\begin{equation*}
A=\left( AA^{T}\right) ^{1/2}O
\end{equation*}%
for some $O\in SO\left( 2\right) $. Hence $E=\left( AA^{T}\right)
^{1/2}S^{1} $ and so, by (\ref{ellipse}), the curve $\gamma $ in $\mathcal{S}%
_{+}$ associated with $E$ is $\gamma =AA^{T}=\alpha \alpha ^{T}+\alpha
^{\prime }\left( \alpha ^{\prime }\right) ^{T}$. Therefore, 
\begin{equation*}
\gamma ^{\prime }=\alpha ^{\prime }\alpha ^{T}+\alpha \left( \alpha ^{\prime
}\right) ^{T}+\alpha ^{\prime \prime }\left( \alpha ^{\prime }\right)
^{T}+\alpha ^{\prime }\left( \alpha ^{\prime \prime }\right) ^{T}\text{.}
\end{equation*}%
Evaluating at $t=0$ one has, using (\ref{normaliz}) and (\ref{cacurvature}),
that%
\begin{equation}
\gamma _{0}^{\prime }=\left( 
\begin{array}{cc}
0 & 0 \\ 
0 & \varkappa _{0}%
\end{array}%
\right) \text{.}  \label{gamma0prima}
\end{equation}%
In particular, $\left\Vert \gamma _{0}^{\prime }\right\Vert =-\det \left(
\gamma _{0}^{\prime }\right) =0$. Since $t_{o}$ was arbitrary, $E$ is a null
curve in $\mathcal{E}$.

In order to prove the second assertion we compute 
\begin{equation*}
\gamma ^{\prime \prime }=\alpha ^{\prime \prime }\alpha ^{T}+2\alpha
^{\prime }\left( \alpha ^{\prime }\right) ^{T}+\alpha \left( \alpha ^{\prime
\prime }\right) ^{T}+\alpha ^{\prime \prime \prime }\left( \alpha ^{\prime
}\right) ^{T}+2\alpha ^{\prime \prime }\left( \alpha ^{\prime \prime
}\right) ^{T}+\alpha ^{\prime }\left( \alpha ^{\prime \prime \prime }\right)
^{T}\text{.}
\end{equation*}%
Now, (\ref{cacurvature}) yields $\alpha ^{\prime \prime \prime }=-\frac{1}{2}%
\varkappa \alpha +\left( \frac{1}{4}\varkappa ^{2}+\frac{1}{2}\varkappa
^{\prime }-1\right) \alpha ^{\prime }$. From this and (\ref{normaliz}) we
have that the first component of $\alpha _{0}^{\prime \prime \prime }$
equals $-\frac{1}{2}\varkappa _{0}$ and 
\begin{equation*}
\gamma _{0}^{\prime \prime }=\left( 
\begin{array}{cc}
0 & -\varkappa _{0} \\ 
-\varkappa _{0} & x%
\end{array}%
\right)
\end{equation*}%
for some number $x$. Since $E\left( t\right) $ and $\gamma _{t}$ correspond
under the isometry (\ref{isometry}), we have by Lemma \ref{DbarD} that 
\begin{equation*}
\left\Vert \left. \frac{DE^{\prime }}{dt}\right\vert _{0}\right\Vert
=\left\Vert \gamma _{0}^{\prime \prime }\right\Vert =-\det \left( \gamma
_{0}^{\prime \prime }\right) =\varkappa _{0}^{2}\text{,}
\end{equation*}%
as desired. The last assertion follows from (\ref{gamma0prima}) and the
definition of the temporal orientation on $\mathcal{S}_{+}$. \hfill $\square 
$

\bigskip

\begin{lemma}
\label{lema}Let $A\in \mathcal{S}_{+}$ and let $v,w$ two vectors in $\mathbb{%
R}^{2}$ such that $v\wedge w>0$, $\left\langle Av,v\right\rangle =1$ and $%
\left\langle Av,w\right\rangle =0$, and let $L=\sqrt{\left\langle
Aw,w\right\rangle }$. Then the curve $\varepsilon :\mathbb{R}\rightarrow 
\mathbb{R}^{2}$ defined by 
\begin{equation}
\varepsilon \left( s\right) =\cos s~Av+\sin s~Aw/L  \label{epsilon}
\end{equation}%
parametrizes the ellipse $E_{A}$ counterclockwise and its curvature at $s=0$
is 
\begin{equation}
\kappa \left( 0\right) =\frac{v\wedge w}{\left( v\wedge Aw\right) \left\vert
v\right\vert }\text{.}  \label{curvature}
\end{equation}
\end{lemma}

\noindent \textbf{Proof. }We verify that $\varepsilon \left( s\right) \in
E_{A}$ for all $s$. Indeed, 
\begin{eqnarray*}
\left\langle A^{-1}\varepsilon \left( s\right) ,\varepsilon \left( s\right)
\right\rangle &=&\left\langle \cos s~v+\sin s~w/L,\cos s~Av+\sin
s~Aw/L\right\rangle \\
&=&\cos ^{2}s~\left\langle v,Av\right\rangle +\sin ^{2}s~\left\langle
Aw,w\right\rangle /L^{2}+2\sin s\cos s~\left\langle w,Av\right\rangle \\
&=&\cos ^{2}s+\sin ^{2}s=1\text{.}
\end{eqnarray*}%
We have that $\varepsilon \wedge \varepsilon ^{\prime }=\left( \det A\right)
v\wedge w/L$ and so $\varepsilon $ parametrizes $E_{A}$ counterclockwise. We
compute 
\begin{equation*}
\kappa \left( 0\right) =\frac{\varepsilon ^{\prime }\wedge \varepsilon
^{\prime \prime }}{\left\vert \varepsilon ^{\prime }\right\vert ^{3}}\left(
0\right) =\frac{Aw/L\wedge \left( -Av\right) }{\left( \left\vert
Aw\right\vert /L\right) ^{3}}=\frac{\left( Av\wedge Aw\right) \left\langle
Aw,w\right\rangle }{\left\vert Aw\right\vert ^{3}}\text{.}
\end{equation*}%
Since $A$ is symmetric we have that $\left\langle Aw,v\right\rangle =0$ and
we may suppose that $Aw=i\mu v$ for some $\mu \in \mathbb{R}$. Here $i$
denotes counterclockwise rotation through a right angle. Therefore 
\begin{equation*}
\kappa \left( 0\right) =\frac{\left( Av\wedge i\mu v\right) \left\langle
i\mu v,w\right\rangle }{\left\vert i\mu v\right\vert ^{3}}=\frac{%
\left\langle v,Av\right\rangle \left( v\wedge w\right) }{\mu \left\vert
v\right\vert ^{3}}=\frac{v\wedge w}{\left\vert v\right\vert ^{3}}
\end{equation*}%
and (\ref{curvature}) holds since $v\wedge Aw=v\wedge i\mu v=\mu \left\vert
v\right\vert ^{2}$. \hfill $\square $

\bigskip

The following theorem is a partial converse of Theorem \ref{TeoP}.

\begin{theorem}
Let $\gamma :\left( a,b\right) \rightarrow \mathcal{S}_{+}$ a regular null
curve with spatial acceleration, that is, $\left\Vert \frac{D\gamma ^{\prime
}}{dt}\right\Vert >0$, and let $E_{t}$ be the ellipse associated with $%
\gamma \left( t\right) $ via the bijection \emph{(\ref{isometry})}. Then
there exists a $0$-convex curve $\alpha :\left( a,b\right) \rightarrow 
\mathbb{R}^{2}$ such that either $t\mapsto E_{t}$ or $t\mapsto E_{-t}$
osculates $\alpha $ at $t$, for all $t\in \left( a,b\right) $.
\end{theorem}

\noindent \textbf{Proof.} We have $\left\Vert \gamma ^{\prime }\right\Vert
=-\det \gamma ^{\prime }=0$. Since $\gamma ^{\prime }\neq 0$, for each $t$,
the kernel of $\gamma ^{\prime }\left( t\right) $ has dimension one. Let $v$
be a smooth curve such that $\gamma ^{\prime }\left( t\right) v\left(
t\right) =0$ for all $t$, normalized in such a way that $\left\langle \gamma
v,v\right\rangle =1$, and let $\alpha \left( t\right) =\gamma _{t}v_{t}$. In
particular, $\alpha ^{\prime }=\gamma ^{\prime }v+\gamma v^{\prime }=\gamma
v^{\prime }$.

By differentiating $1=\left\langle \gamma v,v\right\rangle $ and using that $%
\gamma $ is symmetric we have 
\begin{equation}
0=\left\langle \gamma v^{\prime }+\gamma ^{\prime }v,v\right\rangle
+\left\langle \gamma v,v^{\prime }\right\rangle =\left\langle \gamma
v^{\prime },v\right\rangle +\left\langle \gamma v,v^{\prime }\right\rangle
=2\left\langle \gamma v,v^{\prime }\right\rangle \text{.}  \label{tal}
\end{equation}

Now we verify that $\alpha \ $is 0-convex. Since $\gamma ^{\prime }\neq 0$
is symmetric and singular, and $\gamma ^{\prime }v=0$, there exist a nowhere
vanishing smooth function $\lambda :\left( a,b\right) \rightarrow \mathbb{R}$
such that 
\begin{equation*}
\gamma ^{\prime }=\lambda uu^{T}\text{,}
\end{equation*}%
where $u=iv$ (a column vector). We compute 
\begin{equation*}
\gamma ^{\prime \prime }=\lambda ^{\prime }uu^{T}+\lambda \left(
uu^{T}\right) ^{\prime }=\frac{\lambda ^{\prime }}{\lambda }\gamma ^{\prime
}+\lambda \left( B+B^{T}\right) \text{,}
\end{equation*}%
where $B=u^{\prime }u^{T}$. On the other hand, $\left\langle \gamma ^{\prime
},\gamma ^{\prime }\right\rangle =0$ implies $\left\langle \gamma ^{\prime
\prime },\gamma ^{\prime }\right\rangle =0$. Hence, 
\begin{equation*}
0=\left\langle \gamma ^{\prime \prime },\gamma ^{\prime }\right\rangle =%
\frac{\lambda ^{\prime }}{\lambda }\left\langle \gamma ^{\prime },\gamma
^{\prime }\right\rangle +\lambda \left\langle B+B^{T},\gamma ^{\prime
}\right\rangle =\lambda \left\langle B+B^{T},\gamma ^{\prime }\right\rangle 
\text{.}
\end{equation*}%
Therefore, by Lemma \ref{DbarD}, 
\begin{equation*}
0<\left\Vert \frac{D}{dt}\gamma ^{\prime }\right\Vert =\left\Vert \gamma
^{\prime \prime }\right\Vert =\lambda ^{2}\left\Vert B+B^{T}\right\Vert
=-\lambda ^{2}\det \left( B+B^{T}\right) =\lambda ^{2}\left( u\wedge
u^{\prime }\right) ^{2}=\lambda ^{2}\left( v\wedge v^{\prime }\right) ^{2}%
\text{.}
\end{equation*}%
Consequently, $v\wedge v^{\prime }\neq 0$. We may suppose that $v\wedge
v^{\prime }>0$ (otherwise, we can substitute $\gamma \left( t\right) $ and $%
v\left( t\right) $ with $\overline{\gamma }\left( t\right) =\gamma \left(
-t\right) $ and $\overline{v}\left( t\right) =v\left( -t\right) $,
respectively, and in this case $t\mapsto \overline{E}\left( t\right)
=E\left( -t\right) $ will be the curve of ellipses we were looking for).

We know that $\alpha ^{\prime }=\gamma v^{\prime }$. Hence, 
\begin{equation*}
\alpha \wedge \alpha ^{\prime }=\gamma v\wedge \gamma v^{\prime }=\left(
\det \gamma \right) \left( v\wedge v^{\prime }\right) >0
\end{equation*}%
and so $\alpha $ is parametrized counterclockwise. By (\ref{tal}), since $%
\gamma $ is symmetric, we have $\alpha ^{\prime }=\ell iv$ for some nowhere
zero function $\ell $. Hence, $\alpha ^{\prime \prime }=\ell ^{\prime
}iv+\ell iv^{\prime }$. Thus, 
\begin{equation*}
\alpha ^{\prime }\wedge \alpha ^{\prime \prime }=\ell iv\wedge \left( \ell
^{\prime }iv+\ell iv^{\prime }\right) =\ell ^{2}\left( iv\wedge iv^{\prime
}\right) =\ell ^{2}\left( v\wedge v^{\prime }\right) >0\text{.}
\end{equation*}

Therefore $\alpha $ is 0-convex. Let $\varepsilon $ be as in (\ref{epsilon}%
), with $A=\gamma \left( t\right) $ and $w=v^{\prime }\left( t\right) $, and
let $E\left( t\right) $ be the ellipse parametrized by $\varepsilon $. Then $%
\varepsilon \left( 0\right) =\alpha \left( t\right) $ and $\varepsilon
^{\prime }\left( 0\right) $ is a positive multiple of $\gamma \left(
t\right) v^{\prime }\left( t\right) =\alpha ^{\prime }\left( t\right) $.
Thus, in order to prove that $E\left( t\right) $ osculates $\alpha $ at $t$
it remains only to see that the curvature of $E_{t}$ at $\alpha \left(
t\right) $ coincides with the curvature of $\alpha $ at $t$. The latter is $%
\alpha ^{\prime }\left( t\right) \wedge \alpha ^{\prime \prime }\left(
t\right) /\left\vert \alpha ^{\prime }\left( t\right) \right\vert ^{3}$. We
compute%
\begin{equation*}
\frac{\alpha ^{\prime }\wedge \alpha ^{\prime \prime }}{\left\vert \alpha
^{\prime }\right\vert ^{3}}=\frac{\ell iv\wedge \left( \ell ^{\prime
}iv+\ell iv^{\prime }\right) }{\left\vert \ell iv\right\vert ^{3}}=\frac{%
iv\wedge iv^{\prime }}{\ell \left\vert v\right\vert ^{3}}=\frac{v\wedge
v^{\prime }}{\left( v\wedge \gamma v^{\prime }\right) \left\vert
v\right\vert }
\end{equation*}%
(the last equality follows since $v\wedge \gamma v^{\prime }=v\wedge \alpha
^{\prime }=v\wedge \ell iv=\ell \left\vert v\right\vert ^{2}$), which
coincides with the curvature of $E_{t}$ at $\alpha \left( t\right) $ by
Lemma \ref{lema}. Consequently $E\left( t\right) $ is the osculating ellipse
to $\alpha $ at $t$. \hfill $\square $

\bigskip

The following theorem is the centro-affine analogue of Theorems 5.2 and 7.3
in \cite{LO'H}. In the context of the last paragraph of the introduction, we
consider as $\mathcal{P}_{o}$ the set of all 0-convex paths with nowhere
vanishing centro-affine curvature and give the centro-affine arc length
element as the $\frac{1}{2}$-dimensional length of the curve of osculating
ellipses.

\begin{theorem}
\label{caelement}Let $c$ be a $0$-convex path in $\mathbb{R}^{2}$ with
nowhere vanishing centro-affine curvature. Let $\alpha :I\rightarrow \mathbb{%
R}^{2}$ be any parametrization of $c$ \emph{(}not necessarily standard
centro-affine\emph{)}. For each $t\in I$, let $E\left( t\right) $ be the
ellipse centered at the origin osculating $\alpha $ at $t$. Then the null
curve $E$ in $\mathcal{E}$ has spatial acceleration, that is $\left\Vert 
\frac{D}{dt}E^{\prime }\right\Vert >0$ and the $1$-form $\tau _{c}$ on $c$
is well defined by 
\begin{equation*}
\alpha ^{\ast }\tau _{c}=\left\Vert \frac{DE^{\prime }}{dt}\right\Vert
^{1/4}dt\text{.}
\end{equation*}%
Moreover, the map $c\mapsto \tau _{c}$ is invariant under the action of $G$.
\end{theorem}

\noindent \textbf{Proof. }Notice that the curve $E$ is null by the first
assertion of Theorem \ref{TeoP}. Suppose that $\beta \left( s\right) =\alpha
\left( \phi \left( s\right) \right) $ is a reparametrization of $\alpha $
and let $F\left( s\right) =E\left( \phi \left( s\right) \right) $ be the
osculating ellipse of $\beta $ at $s$. We compute $F^{\prime }\left(
s\right) =E^{\prime }\left( \phi \left( s\right) \right) \phi ^{\prime
}\left( s\right) $ and%
\begin{equation*}
\frac{DF^{\prime }}{ds}=\frac{DE^{\prime }}{dt}\left( \phi \right) \left(
\phi ^{\prime }\right) ^{2}+E^{\prime }\left( \phi \right) \phi ^{\prime
\prime }\text{.}
\end{equation*}%
Now $\left\langle E^{\prime },E^{\prime }\right\rangle =0$ and this implies $%
\left\langle \frac{DE^{\prime }}{dt},E^{\prime }\right\rangle =0$. Hence 
\begin{equation}
\left\Vert \frac{DF^{\prime }}{ds}\right\Vert =\left\Vert \frac{DE^{\prime }%
}{dt}\left( \phi \right) \right\Vert \left( \phi ^{\prime }\right) ^{4}\text{%
.}  \label{DFDE}
\end{equation}%
If $\beta $ is a standard centro-affine reparametrization of $\alpha $, we
have by hypothesis and the second assertion of Theorem \ref{TeoP} that $%
\left\Vert \frac{DF^{\prime }}{ds}\right\Vert $ is positive, and by (\ref%
{DFDE}) $\left\Vert \frac{DE^{\prime }}{ds}\right\Vert $ is also so.
Therefore%
\begin{equation*}
\left\Vert \frac{DF^{\prime }}{ds}\left( s\right) \right\Vert
^{1/4}ds=\left\Vert \frac{DE^{\prime }}{dt}\left( \phi \left( s\right)
\right) \right\Vert ^{1/4}\phi ^{\prime }\left( s\right) ~ds\text{,}
\end{equation*}%
and this implies that the $1$-form $\tau _{c}$ on $c$ is well defined, since 
\begin{equation*}
\beta ^{\ast }\left( \alpha ^{-1}\right) ^{\ast }=\left( \alpha ^{-1}\beta
\right) ^{\ast }=\phi ^{\ast }\text{\ \ \ \ \ \ \ and\ \ \ \ \ \ \ \ }\phi
^{\ast }dt=\phi ^{\prime }\left( s\right) ~ds\text{.}
\end{equation*}

Finally, we show that $\tau _{c}$ is $G$-invariant. We have to check that $%
\tau _{c}=A^{\ast }\tau _{Ac}$ for any $A\in G$, or equivalently, that $%
\alpha ^{\ast }\tau _{c}=\left( A\alpha \right) ^{\ast }\tau _{Ac}$ for a
parametrization $\alpha $ of $c$, that is, 
\begin{equation*}
\left\Vert \frac{DE^{\prime }}{dt}\right\Vert ^{1/4}dt=\left\Vert \frac{%
D\left( AE\right) ^{\prime }}{dt}\right\Vert ^{1/4}dt\text{,}
\end{equation*}%
($AE\left( t\right) $ osculates $A\alpha $ at $t$) and this is true since $A$
acts by isometries on $\mathcal{E}$. \hfill $\square $

\bigskip

\noindent Marcos Salvai

\noindent \textsc{ciem - f}a\textsc{maf}

\noindent Conicet - Universidad Nacional de C\'{o}rdoba

\noindent Ciudad Universitaria, 5000 C\'{o}rdoba, Argentina

\noindent salvai@famaf.unc.edu.ar

\end{document}